\documentclass[12pt,a4paper]{article}

\usepackage{setspace}
\usepackage{authblk}

\usepackage{float}
\restylefloat{table}

\usepackage{graphicx}	 			
\usepackage{amsmath}
\usepackage{amssymb}
\usepackage{mathrsfs}
\usepackage{amsthm}
\usepackage{graphicx}	

\usepackage{cancel}

\usepackage{epsfig,amstext,amsfonts,color}

\usepackage{latexsym}
\usepackage{hyperref}[bookmarksopen]

\renewcommand{\mathbf}{\boldsymbol}

\date{}

\makeatletter
\renewcommand\@biblabel[1]{}
\makeatother

\begin{document}

\begin{center}
\textbf{\Large A multivariate generalization of Hall's Theorem for Edgeworth expansions of bootstrap distributions}\\
\vskip 0.1truein
\textbf{\large Andrew T.A. Wood}\\
\textbf{\large (Australian National University)}
\end{center}
\vskip 0.2truein
\begin{center}
\textbf{Abstract}
\end{center}
Theorem 5.1 in the monograph by  Hall (1992) provides rigorous in-probability justification of Edgeworth expansions of bootstrap distributions.  Proving this result was rather challenging because bootstrap distributions do not satisfy the classical Cram\'er  condition and therefore classical methods for justifying Edgeworth expansions, e.g. Bhattacharya and Rao (1976) and Bhattacharya and Ghosh (1978), are not available.  Hall's (1992) theorem is for a univariate statistic which can be expressed as a smooth function of means, though the underlying population can be multivariate.  However, there are a number of applications where a multivariate version of Hall's theorem is needed, and generalizing the proof from the univariate case to the multivariate case is not immediate.  Our primary purpose in this article is to fill this gap by stating a multivariate version of the theorem and sketching the modifications to the proof of Hall's (1992) Theorem 5.1 that are needed.


\vskip 0.2truein


\vfill \eject

\section{Introduction}
The classical sufficient conditions for the validity of an Edgeworth expansion of the sample mean of an IID sequence of random vectors, in the sense that the actual error in the Edgeworth approximation is the same as the nominal error (see below for specifics), are:  (i) sufficiently many moments exist and (ii) Cram\'er's condition holds (see (\ref{cramer}) below). An important technical point is that the bootstrap distribution does not satisfy Cram\'er's condition (ii), due to its discreteness.   However, Cram\'er's condition is, in fact, far stronger than is necessary for an Edgeworth expansion to be valid.  To gain a sense of how far from necessary Cram\'er's condition is for rigorous Edgeworth expansion validity, see Theorem 2.1 in Booth et al. (1994).  However, when Cram\'er's condition does not hold, proving Edgeworth expansion validity becomes a far more challenging problem.  So far as we are aware, Peter Hall was the only person who ever seriously grappled with this problem in its full generality.

In Theorem 5.1 of Hall's (1992) monograph, a proof of the validity of an Edgeworth approximation of arbitrary order is given for the bootstrap distribution of a univariate statistic that may be expressed as a smooth function of means.  The proof of Hall's Theorem 5.1 is long and very technical; the somewhat abbreviated proof of Theorem 5.1 given in Hall (1992), which omits many details, is 21 pages.  It takes as a starting point the  detailed analyses of Bhattacharya and Rao (1976) and Bhattacharya and Ghosh (1978) of the remainder terms which arise in Edgeworth expansions in classical settings.  Hall's (1992) immense theoretical contribution is to prove that, under sufficient moment conditions, under sufficient smoothness of the function in the smooth function model, and assuming the underlying population satisfies Cram\'er's condition, the following holds: with probably approaching 1 at a suitably fast rate, the order of the relevant remainder term in the Edgeworth expansion is actually equal to the nominal (heuristic) order.  In this sense, an Edgeworth approximation of arbitrary order  is valid for the bootstrap distribution in probability, under suitably strong versions of the assumptions mentioned above.

The reason why Hall's (1992) Theorem 5.1 is a theoretical result with practical relevance is that it enables one to make rigorously justified statements about the higher-order accuracy of bootstrap procedures; e.g.  statements such as ``the coverage error of a certain type of bootstrap confidence region is $O(n^{-2})$'' or ``the $p$-value for a certain type of bootstrap hypothesis test is accurate to order $O_p(n^{-3/2})$''.

As already mentioned, Hall's (1992) Theorem 5.1 is for a univariate statistic that can be expressed as a smooth function of vector means.  However, there are many  examples of interest where the statistic expressible as a smooth function of vector means is multivariate, not univariate.  Examples of such statistics are given in Fisher et al. (1996), Amaral et al. (2007)  and Hingee et al. (2026). The  purpose of this article is to provide a multivariate generalization of Hall's (1992) Theorem 5.1 which justifies the claims in these papers. 

A complete, self-contained proof of our generalization of Hall's (1992) Theorem 5.1 would be too long and would involve too much repetition of Hall's arguments to be a worthwhile undertaking.  For these reasons our strategy is as follows.  After introducing notation and assumption and then stating the main result in Section 2, we then explain in Section 3 how to generalize the definitions of the main ingredients used in  Hall's proof.  Some of these definitions go through without change but other definitions require modification.  Then, in Section 4, we state and prove a new result, Proposition 1, which generalizes Hall's (1992) Lemma 5.3 to the multivariate case.  Formulating and proving our Proposition 1  actually turns out to be the biggest obstacle in generalizing Hall's (1992) proof  of Theorem 5.1.  In Section 5 we discuss certain sets which arise in the proof of Hall's Theorem 5.1.  In Section 6, we give a brief overview of the proof and then conclude with discussion in Section 7.

\section{Notation, assumptions  and main result}

\subsection{Notation and assumptions}

First, let us consider a class of Borel subsets, $\mathcal{B}$, of $\mathbb{R}^q$ with the following properties: (i) for $x, y \in \mathbb{R}^q$ and any $B \in \mathcal{B}$,
\begin{equation}
   y+B=\{x+y: \, x \in B\} \in \mathcal{B}, 
\label{curlyB_1}
\end{equation}
i.e. $\mathcal{B}$ is closed under translations; and (ii)
denoting the multivariate normal probability density function with mean $0$ ($q \times 1$)  and covariance $V$ ($q \times q$) by $\phi_{0,V}$, we have, for some $C>0$, 
\begin{equation}
\sup_{\lambda_{\max}(V) \leq C, \, \lambda_{\min}(V)> C^{-1}}\, \sup_{B \in \mathcal{B}} \int_{(\partial B)^\epsilon}
 \phi_{0,V}(x) dx = O(\epsilon), \hskip 0.2truein \epsilon \downarrow 0,
 \label{curlyB_2}
\end{equation}
where $\partial B$ is the boundary of $B$, 
\[
(\partial B)^\epsilon = \left \{ y \in \mathbb{R}^q: \vert \vert x-y \vert \vert <\epsilon \hskip 0.1truein \textrm{for some} \hskip 0.1truein x \in \partial B \right \}
\]
and $\lambda_{\max}(V)$ and $\lambda_{\min}(V)$ are, respectively, the maximum and minimum eigenvalues of a symmetric matrix $V$.

We now  specify Cram\'er's condition.  If $\chi(t)=E[\exp(i t^\top X)]$ is the characteristic function of a random vector $X \sim F$, then $X$ (and the probability distribution $F$) are said to satisfy Cram\'er's condition if
\begin{equation}
\sup_{\vert \vert t \vert \vert \rightarrow \infty} \vert \chi(t) \vert <1.
\label{cramer}
\end{equation}

Assumptions (\ref{curlyB_1}), (\ref{curlyB_2}) and (\ref{cramer}) are all rather mild.  If necessary, we can always enlarge the set of subsets $\mathcal{B}$ so that (\ref{curlyB_1}) holds.  In equation (3) on page 2 of Bharracharya and Rao (1976), it is stated that (\ref{curlyB_2}) holds if $V=I_q$ and $\mathcal{B}$ is the set of convex sets in $\mathbb{R}^q$; for many applications of Theorem 1 below, an appropriate  choice of $\mathcal{B}$ will be the subset of the balls in a Euclidean space.  Finally, a sufficient condition for Cram\'er's condition to hold is that, for some integer $n_0$ satisfying $1 \leq n_0 < \infty$, a sum of $n_0$ copies of $X$ has a component which is absolutely continuous with respect to Lebesgue measure on $\mathbb{R}^q$.

In some of what follows we will consider $k$-sample problems in which there are $k$ population distributions of interest defined on $\mathbb{R}^d$.  Let $\mathcal{P}_j$ denote population distribution $j$, $j=1, \ldots , k$, and consider samples 
$\mathcal{X}_j=\{X_{j1}, \ldots , X_{jn_j}\}$ drawn randomly from population $\mathcal{P}_j$, and write $\bar{X}_j$ for the sample mean based on $\mathcal{X}_j$.  Let $\mathcal{X}_j^\ast$ denote a bootstrap sample with sample size $n_j$ drawn randomly (i.e. with equal probabilities of selection) with replacement from $\mathcal{X}_j$.
Let us write $\bar{X}_1^\ast, \ldots , \bar{X}_k^\ast$ for the $k$ sample mean
vectors based, respectively, on the bootstrap samples $\mathcal{X}_1^\ast, \ldots , \mathcal{X}_k^\ast$. 


In order to avoid being distracted by uninteresting degeneracies, we assume that $n_j$, $j=1, \ldots , k$, satisfy 
\begin{equation}
\limsup_{n \rightarrow \infty}\frac{\max(n_1, \ldots, n_k)}{\min(n_1, \ldots , n_k) } < \infty,
\end{equation}
where $n=(n_1+ \cdots + n_k)$.

In what follows we require a more general definition of the non-random function $A(x)=\{g(x)-g(\mu)\}/h(x)$ and random function $\hat{A}(x) =\{g(x)-g(\bar{X})\}/h(x)\}$ defined on page 238 of Hall (1992), where $\mu=E(X)$ is the mean of a single observation vector and both $A$ and $\hat{A}$ are real-valued.  For the more general setting of Theorem 1, we define $A$ and $\hat{A}$ as follows:
\begin{equation}
A=g(x_1, \ldots , x_k; \mu_1, \ldots , \mu_k) \hskip 0.2truein
\textrm{and} \hskip 0.2truein \hat{A}=g(x_1, \ldots , x_k; \bar{X}_1, \ldots , \bar{X}_k),
\label{new_A}
\end{equation}
where $x_j \in \mathbb{R}^d$, $\bar{X}_j$ is the sample mean vector based on sample $\mathcal{X}_j$ drawn randomly from infinite population $\mathcal{P}_j$ and $\mu_j = E(X_j)$ where $X_j \sim \mathcal{P}_j$, $j=1, \ldots , k$.  Note that the range of $A$ and $\hat{A}$ lies within $\mathbb{R}^q$.

Finally, we recall the functions and signed measures defined in pages 51-57 of Bhattacharya and Rao (1976):
\begin{equation}
P_j(-\phi_{0,V}: \{\kappa_{r}\})(x), \hskip 0.2truein j=0,1,\ldots
\label{P_phi}
\end{equation}
and the signed measure 
\begin{equation}
P_j(-\Phi_{0,V}: \{ \kappa_r \})(\mathcal{R}) = \int_{\mathcal{R}}
P_j(-\phi_{0,V}: \{\kappa_r \})(x)dx, \hskip 0.2truein j=0,1,\ldots
\label{P_Phi}
\end{equation}
where $\phi_{0,V}(x)$ is the $d$-dimensional multivariate normal density with mean vector $0$ and $d \times d$ positive definite covariance matrix $V$, $\Phi_{0,V}$ is the Gaussian corresponding probability measure on $\mathbb{R}^d$ with mean vector $0$ and covariance matrix $V$, $\{\kappa_r  \}$ is a cumulant sequence, and $\mathcal{R} \subseteq \mathbb{R}^d$ is any Lebesgue measurable set. The precise definition of (\ref{P_phi}) is complex and is not required here.  Key properties of (\ref{P_phi}) and (\ref{P_Phi}) are derived on pages 51-57 of Bhattacharya and Rao (1976).  Several points are worth noting.  First, define
\begin{equation}
\check{P}_j(x)=\check{P}_j(0,V, \{\kappa_r\})(x)=\frac{P_j(-\phi_{0,V}: \{\kappa_r\})(x)}{\phi_{0,V}(x)}.
\label{check_P_j}
\end{equation}
Then $\check{P}_j(0,V,\{\kappa_r\})(x)$ is a polynomial in the components of $x=(x_1, \ldots , x_d)^\top$.  The polynomials $\check{P}_j(x)$ are of degree $3j$ and are even functions of $x$ for even $j$ and odd functions of $x$ for odd $j$ in the sense that
$\check{P}_j(-x) =(-1)^j \check{P}_j(x)$.
Moreover, the coefficients of $\check{P}_j$ only involves joint cumulants (equivalently, joint moments) up to order $j+2$.

\subsection{The main result}

We now state the main result of this article.

\noindent \textbf{Theorem 1}.  \textit{Let $\lambda>0$ be given.  Suppose that, for each of the $k$ populations, a random observation from population $j$, $X_j$ say, ($d \times 1$) is such that  $E[\vert \vert X_j \vert \vert^s]<\infty$ where $s=2 \lambda  \max(2 \nu+3, d+1)$, where $\nu\geq 1$ is a positive integer whose role is indicated below, and assume that each $X_j$ satisfies Cram\'er's condition (\ref{cramer}).  Let $A: \mathbb{R}^{dk} \rightarrow \mathbb{R}^q$ defined in (\ref{new_A}) denote a  function all of whose partial derivatives of order $\nu+3$ or less are continuous in a neighbourhood of $\mu=(\mu_1^\top , \ldots , \mu_k^\top)^\top$  and let $\mathcal{B}$ denote a subset of the Borel subsets of $\mathbb{R}^q$ such that (\ref{curlyB_1}) and (\ref{curlyB_2}) are satisfied.  Then there exists a constant $C_0>0$ such that
\begin{align}
P\bigg [ \sup_{B \in \mathcal{B}}\bigg \vert P&\left \{n^{1/2} \hat{A}(\bar{X}_1^\ast, \ldots , \bar{X}_k^\ast) \in B \vert \mathcal{X} \right \} \nonumber \\
&-\sum_{j=0}^\nu
n^{-j/2} P_j(-\Phi_{0,I_q}: \{\hat{\eta_r}\})(B) \bigg \vert >
C_0n^{-(\nu+1)/2} \bigg ] = O(n^{-\lambda}),
\label{general_hall}
\end{align}
where $\{\hat{\eta}_r\}$ is actually an approximate cumulant sequence for $n^{1/2}A(\bar{X}_1^\ast, \ldots , \bar{X}_k^\ast)$, and
\[
P \left \{ \max_{1 \leq j \leq \nu} \, \max_{x \in \mathbb{R}^{dk}}
\left ( 1+\vert \vert x \vert \vert^{-3j} \right) \vert \check{P}_j(x) \vert  >C_0 \right \}=O(n^{-\lambda}),
\]
where $\hat{A}$ is defined in (\ref{new_A}).}


  In (\ref{general_hall}), $\{\hat{\eta}_r\}$ is the cumulant sequence for a sufficiently high-order polynomial approximation in several variables to $n^{1/2} A(\bar{X}_1^\ast, \ldots , \bar{X}_k^\ast)$.  The relevant approximate cumulants have expansions in powers of $n^{-1}$; see Theorem 2.1 in Hall (1992) for key insights and relevant discussion.  

Theorem 1 is a generalization of Theorem 5.1 in Hall (1992) in two senses: Theorem 5.1 in Hall (1992)  corresponds to $k=1$ and to $q=1$, with 
$\mathcal{B}=\left \{(-\infty , x]: \, x \in \mathbb{R}  \right \}$.  In Theorem 1 above we allow $k \geq 1$ and/or $q \geq 1$.  The generalization which requires new ingredients in the proof is $q>1$.

The assumption $\mathbb{E}[\vert \vert X_jvert \vert^s ]$ is very likely to be far stronger than necessary; see Point 4 is the discussion in Section 7.


Before sketching the proof of Theorem 1 in subsection \ref{Proof_of_S1} we present some modified definitions from Hall (1992) and supporting results in the following subsections, including a statement and proof of Proposition 1 in Section 4, a new result which is needed in the proof of Theorem 1.

\section{Modified definitions in Hall's Theorem 5.1 \label{modsABC}}

As mentioned above, Theorem 1 is a generalization of Theorem 5.1 in Hall (1992).  As a consequence, certain terms that arise in the proof need to be redefined.  We present the modified definitions in this section.

\subsection{Definition of the sets $\mathcal{R}(B)$, $B \in \mathcal{B}$. \label{curlyR_B}}
In the proof of Theorem 5.1 in Hall (1992), a central role is played by the sets
\begin{equation}
\mathcal{R}(t)=\{x \in \mathbb{R}^d: n^{1/2}\hat{A}(\bar{X}+n^{-1/2}x) \leq t\}, \hskip 0.2truein t \in (- \infty , \infty).
\label{Hall77}
\end{equation}
In the proof of Theorem 1, we need to work with a more general class of sets defined as follows:
\begin{equation}
    \mathcal{R}(B)=\{x=(x_1^\top, \ldots , x_k^\top)^\top \in \mathbb{R}^{dk}: n^{1/2}\hat{A}(\bar{X}_1+n^{-1/2}x_1, \ldots , \bar{X}_k + n^{-1/2}x_k) \in B\},
    \label{newguy}
\end{equation}
where $B \in \mathcal{B}$ and each $x_j \in \mathbb{R}^d$, $j=1, \ldots , k$.  Note that (\ref{newguy}) generalizes (\ref{Hall77}) in two senses.  First, the domain of $\hat{A}$ in (\ref{Hall77}) lies within $\mathbb{R}^d$ whereas in (\ref{newguy}) the domain lies within $\mathbb{R}^{dk}$.  Second, the range of $\hat{A}$ in  (\ref{Hall77}) is $\mathbb{R}$ whereas the range of $\hat{A}$ in (\ref{newguy}) is $\mathbb{R}^q$.

The definition of the set $\mathcal{R}^\dagger$ given at the top of page 247 of Hall (1992) also needs modification.  The appropriate definition of $\mathcal{R}^\dagger$ for the proof of Theorem 1 is:
\begin{align}
\mathcal{R}^\dagger (B)&=\bigg \{\left (  (\hat{V}_1^{-1/2}x_1 - n^{1/2}E(Y_1 \vert \mathcal{X}_1))^\top, \ldots ,  (\hat{V}_k^{-1/2} x_k- n^{1/2} E(Y_k\vert \mathcal{X}_k))^\top \right )^\top \in \mathbb{R}^{kd} \nonumber \\
&\hskip 1.0truein : x=(x_1^\top , \ldots , x_k^\top )^\top \in \mathcal{R}(B) \bigg \},
\label{curlyR^dagger}
\end{align}
for $B \in \mathcal{B}$, where $\hat{V}_j$ is the sample covariance matrix for sample $j$, denoted $\mathcal{X}_j$, $j=1, \ldots , k$.  The presence of $n^{1/2}$ does not cause a problem in (\ref{curlyR^dagger}) as $E(Y_{ji}\vert \mathcal{X}_j)$ is small due to the truncation in the definition of $Y_{ji}$ in (\ref{eq:Yji}) below. 

\subsection{Definition of the signed measure $R$ and sets $\mathcal{R}(B)$ and $\mathcal{R}^\dagger(B)$} 

It is  necessary to redefine the signed measure $R$ defined in line 5 on page 247 of Hall (1992). In Hall's proof, $R$ is a signed measure on $\mathbb{R}^d$ whereas in our proof the analogue of $R$ is a signed measure on $\mathbb{R}^{kd}$, due to there being $k$ samples.  

Starting with the definitions for a single sample on page 246 of Hall (1992), we modify the definitions as follows.
\begin{equation}
\hat{V}=\textrm{diag} \left \{\hat{V}_1, \ldots , \hat{V}_k \right \} \hskip 0.2truein \textrm{and} \hskip 0.2truein
\hat{V}^\dagger = \textrm{diag} \left \{\hat{V}_1^\dagger , \ldots , \hat{V}_k^\dagger \right \},
\label{hatVhatV^dagger}
\end{equation}
where $\hat{V}_j$ is the sample covariance matrix of sample $\mathcal{X}_j$ and $\hat{V}_j^\dagger$ is the sample covariance matrix of $X_{ji}^\dagger$, $i=1, \ldots , n_j$, with
\begin{equation}
X_{ji}^\dagger = Y_{ji}-E(Y_{ji}\vert \mathcal{X}_j),
\label{X^dagger}
\end{equation}
and $Y_{ji}$ is defined as
\begin{equation}
Y_{ji}=\begin{cases}\hat{V}_j^{-1/2}(X_{ji}^\ast - \bar{X}_j) & \textrm{if} \, \vert \vert \hat{V}^{-1/2}(X_{ji}^\ast - \bar{X}_j) \vert \vert >n^{1/2}\\
0 & \textrm{otherwise},
\end{cases}
\label{eq:Yji}
\end{equation}
where the $X_{ji}^\ast$ are resampled with equal weights (i.e. resampled randomly with replacement) from sample $\mathcal{X}_j$.  

The generalized definitions of the discrete probability measures $Q$ and $Q^\dagger$ are now given.  In the setting of Theorem 1, $Q$ and $Q^\dagger$ are product measures.  Consider sets $S_1, \ldots , S_k \subset \mathbb{R}^d$ and define
$S = S_1 \times \cdots S_k \subset \mathbb{R}^{kd}$.  $S$ is known as a cylinder set.  On the class of cylinder sets we define the product measures
\begin{equation}
Q(S)=\prod_{j=1}^k Q_j(S_j) \hskip 0.2truein \textrm{and} \hskip 0.2truein Q^\dagger (S) =\prod_{j=1}^k Q_j^\dagger (S_j),
\label{bigQ}
\end{equation}
where each $Q_j$ and $Q_j^\dagger$ is based on sample $\mathcal{X}_j$ and is defined in exactly the same way as Hall (1992, p. 246) defines $Q$ and $Q^\dagger$ for a single sample.  More specifically, $Q_j$ and $Q_j^\dagger$ are  the discrete probability measures generated by 
\[
n_j^{-1/2} \sum_{i=1}^{n_j} (X_{ji}^\ast - \bar{X}_j) \hskip 0.2truein
\textrm{and} \hskip 0.2truein n_j^{-1/2}\sum_{i=1}^{n_j} X_{ji}^\dagger,
\]
where $X_{ji}^\dagger$ is defined in (\ref{X^dagger}).

One technical point: note that the probability measures $Q$ and $Q^\dagger$ are fully determined by the values of the measure on the class of cylinder sets, which is a much smaller class of sets than the Borel sets.

The signed measure given in line 5 on page 247 of Hall (1992) generalizes to
\[
R=Q^\dagger - \sum_{j=1}^{\nu+3} n^{-j/2} P_j(-\Phi_{0,\hat{V}}:
\{ \kappa_{1r}^\dagger \}, \ldots , \{\kappa_{kr}^\dagger\}),
\]
where $Q^\dagger$ is defined by (\ref{bigQ}) and in the discussion below (\ref{bigQ}); $\hat{V}$ is the first matrix defined in (\ref{hatVhatV^dagger}); and $\{\kappa_{jr}^\dagger\}$ is the cumulant sequence for a typical $X_{ji}^\ast$ resampled randomly with equal weights from $\mathcal{X}_j$, the sample from population $j$, $j=1, \ldots , k$.

\subsection{Boundaries and neighborhoods of boundaries}

There are two distinct but related types of boundary that we need to consider in the proof of Theorem 1.  The situation is somewhat simpler in the proof of Theorem 5.1 in Hall (1992) because there $\hat{A}$ is real valued.  The first type is the boundary of the set $B$, denoted $\partial B$, which appears in condition (\ref{curlyB_2}), where $B \subset \mathbb{R}^q$.  The second type of boundary is of the form $\partial \mathcal{R}(B) \subset \mathbb{R}^{kd}$, defined in (\ref{newguy}).  Under mild conditions, in particular that 
\[
f_n(x_1, \ldots , x_k) =n^{1/2} \hat{A}(\bar{X}_1+ n^{-1/2}x_1, \ldots ,\bar{X}_k + n^{-1/2}x_k)
\]
is sufficiently smooth, then
\begin{equation}
\partial \mathcal{R}(B)=\{x =(x_1^\top, \ldots , x_k^\top)^\top
\in \mathbb{R}^{kd}: f_n(x_1, \ldots , x_k) \in \partial B\},
\label{identity}
\end{equation}
for all $B \in \mathcal{B}$.  Moreover, an $\epsilon$-neighborhood of $\partial \mathcal{R}$, denoted $(\partial \mathcal{R}(B))^{\epsilon}$, is defined by
\begin{align*}
(\partial \mathcal{R}(B))^\epsilon
&=\left \{ y \in \mathbb{R}^{kd}: \vert \vert y - x \vert \vert <\epsilon, \, \, x \in \partial \mathcal{R}(B)  \right \}\\
&=\left \{ y \in \mathbb{R}^{kd}: \vert \vert y-x \vert \vert < \epsilon,
\, \, f_n(x_1, \ldots , x_k) \in \partial B\right \}.
\end{align*}

\section{A generalization of Lemma 5.3 in Hall (1992)}

A key role in the proof of Hall's (1992) Theorem 5.1 is played by Hall's (1992) Lemma 5.3.  Proposition 1 below is a multivariate generalization of Hall's (1992) Lemma 5.3 which is needed to prove Theorem 1 above.

\noindent \textbf{Proposition 1}.  \textit{Suppose $Z \sim N_d(0_d, I_d)$ denote a standard $d$-dimensional normal and, for each $\alpha=1, \ldots q$,  consider the polynomials in $x=(x_1, \ldots , x_d)^\top$ given by
\begin{align*}
p_{n\alpha}(x)&=\sum_{i=1}^d c_{\alpha i} x_i +n^{-1/2}
\sum_{i_i=1}^d \sum_{i_2=1}^d c_{\alpha i_1 i_2} x_{i_1} x_{i_2}+ \ldots \\
& \hskip 0.5truein n^{-\nu/2} \sum_{i_1=1}^d \sum_{i_2=1}^d \cdots \sum_{i_\nu=1}^d c_{\alpha i_1 i_2 \cdots i_\nu}x_{i_1} x_{i_2} \cdots x_{i_\nu}.
\end{align*}
whose coefficients $c_{\alpha i}$, $c_{\alpha i_1 i_2}$, \ldots , $c_{\alpha i_1 i_2 \cdots i_\nu}$ satisfy
\[
\max_{\alpha=1, \ldots , q}\, \,  \max_{i_1, i_2, \ldots , i_\nu=1, \ldots , d} \, \, (\vert c_{\alpha i} \vert, \vert c_{\alpha i_1 i_2} \vert, \cdots , \vert c_{\alpha i_1 i_2 \cdots i_\tau} \vert , \cdots , \vert c_{\alpha i_1 i_2, \ldots i_\nu}\vert ) \leq b_1,
\]
for some $0<b_1<\infty$, where $1 \leq \tau \leq \nu$.  Assume that the matrix 
\begin{equation}
W=\left (\sum_{i=1}^d c_{\alpha i } c_{\beta i}\right )_{\alpha, \beta=1}^q=\left ( w_{\alpha \beta} \right )_{\alpha, \beta =1}^q,
\label{full_rank}
\end{equation}
with largest eigenvalue $\lambda_{\max}(W)$ satisfying $\lambda_{\max}(W) \leq b_2$ and smallest eigenvalue $\lambda_{\min}(W)$ satisfying $\lambda_{\min}(W) \geq b_2^{-1}$ for some $b_2 \in (0,\infty)$.
Then, provided $1 \leq q \leq d$, there exists a constant $C>0$ depending only on $\beta$, $b=\max(b_1,b_2)$, $d$ and $q$ such that
\[
\sup_{B \in \mathcal{B}}
P\left \{  \vert \vert Z \vert \vert \leq b\log n \hskip 0.2truein
\textrm{and} \hskip 0.2truein (r_{n1}(Z), \ldots , r_{nq}(Z))^\top \in (\partial B)^{n^{-\beta}}\right \} \leq Cn^{-\beta},
\]
where $C=C(\beta, b, d, q)$.
}

Proposition 1 is a multivariate generalization of Lemma 5.3 of Chapter 5 in Hall (1992).  To see this, take $q=1$, define $B_t=(-\infty , t]$, $\mathcal{B}=\{B_t: t \in \mathbb{R}\}$, and note that the conditions on the coefficients $c_{\alpha i}$ etc.~are essentially the same as in Lemma 5.3 of Hall (1992).

The proof of Proposition 1 makes use of the following lemma which is proved before we start the proof of Proposition 1.

\noindent \textbf{Lemma 1}.  \textit{Let $v_1, \ldots , v_d$ denote $q$-dimensional vectors, define
\[
W_0=\sum_{i=1}^d v_i v_i^\top, \hskip 0.2truein W_j=W_0-v_j v_j^\top, \hskip 0.1truein j=1, \ldots , d,
\]
and suppose $d>q$.  If $W_0$ is positive definite then
\begin{equation}
\max_{j=1, \ldots , d}\lambda_{\min}(W_j) >0,
\label{max_lambda_q}
\end{equation}
where $\lambda_{\min}(W_j)$ is the smallest eigenvalue of the $q \times q$ matrix $W_j$.}

Note that Lemma 1 proves that at least one of the $W_j$ is positive definite.

\noindent \textbf{Proof of Lemma 1}.  If the Lemma is true whenever $d=q+1$ then it must be true for all $d \geq q+1$, so let us focus on the case $d=q+1$. Assume all the $v_j$ are non-zero vectors; if not, just one of the $v_j$ can be the zero vector in which case the conclusion of the lemma holds due to the positive definiteness of $W_0$.  We shall assume that the conclusion of the lemma fails and then derive a contradiction.  Consider $W_1$.  Either (i) it is positive definite or (ii) for $j=2, \ldots , q+1$, $v_j$ is orthogonal to $v_1$, i.e. $v_1^\top v_j=0$; similarly, either (i) $W_2$ is positive definite or (ii) $v_2^\top v_j=0$ for $j\neq 2$.  Repeating this argument, either $W_j$ is positive definite for some $j=1, \ldots , q+1$, in which case the conclusion of the lemma holds or, if not, then $v_j^\top v_k=0$ for all $j\neq k=1, \ldots , q+1$.  Since each $v_j$ is a $q \times 1$ non-zero vector, the latter possibility can not arise, because it is not possible to have $q+1$ non-zero vectors of dimension $q \times 1$ that are mutually orthogonal, so that we have a contraction.  Hence at least one of the $W_j$ must be positive definite and the lemma is proved.  \hfill $\square$

\noindent \textbf{Proof of Proposition 1}.  First, we note that the conclusion of Proposition 1
holds in the case $q = d$.  This follows from the fact that (\ref{full_rank}) has full rank $q$, due to its smallest eigenvalue $\lambda_q(W)$ being positive, and due to the assumptions concerning $\mathcal{B}$, especially (\ref{curlyB_2}). 

For each $q\geq 1$, we use induction on $d$, assuming that $d>q$.  The remainder of the proof is similar to that given in Hall (1986, p.1446); this result is also closely related to Lemma 5.3 in Hall (1992), but Hall's former result is more restricted than the latter.  We prove a more general version of Lemma 5.3 in Hall (1992), in which $q =1$; in our result, the case $q \geq 1$ is included.

Fix $q \geq 1$ and assume now that the result is true for some $d-1\geq q$. Our goal is to prove that it is then true for $d$. We make use of Lemma 1.

 Let $v_j=(c_{1j}, \ldots c_{qj})^\top$, $j=1, \ldots , d$.  Without loss of generality it is arranged that the labeling is chosen such that the $j$ for which $\lambda_q(W_j)$ is maximized is for $j=d$. 
Define $Z^\ast=(Z_1, \ldots , Z_{d-1})^\top$ and $x=(x_1, \ldots , x_d)^\top$ and  $x^\ast=(x_1, \ldots , x_{d-1})^\top$.  Then the polynomial $p_{n \alpha}(x)$ defined below may be re-expressed as follows:
\begin{align*}
p_{n\alpha}(x)&=\sum_{i=1}^d c_{\alpha i} x_i +n^{-1/2}
\sum_{i_i=1}^d \sum_{i_2=1}^d c_{\alpha i_1 i_2} x_{i_1} x_{i_2}+ \ldots \\
& \hskip 0.5truein + n^{-\nu/2} \sum_{i_1=1}^d \sum_{i_2=1}^d \cdots \sum_{i_\nu=1}^d c_{\alpha i_1 i_2 \cdots i_\nu}x_{i_1} x_{i_2} \cdots x_{i_\nu}\\
&=\sum_{i=1}^{d-1} c_{\alpha i}^\ast  x_i +n^{-1/2}
\sum_{i_i=1}^{d-1} \sum_{i_2=1}^{d-1} c_{\alpha ij}^\ast x_{i_1} x_{i_2}+ \ldots \\
& \hskip 0.5truein + n^{-\nu/2} \sum_{i_1=1}^{d-1} \sum_{i_2=1}^{d-1} \cdots \sum_{i_\nu=1}^{d-1} c_{\alpha i_1 i_2 \cdots i_\nu}^\ast x_{i_1} x_{i_2} \cdots x_{i_\nu}\\
&=p_{n\alpha}^\ast (x^\ast \vert x_d),
\end{align*}
where $\alpha = 1, \ldots , q$, $c_{\alpha i}^\ast = c_{\alpha i}^\ast (x_d), \ldots , c_{\alpha i_1 \cdots i_\nu}^\ast=c_{\alpha i_1 \cdots i_\nu}^\ast (x_d)$ are all polynomials in $x_d$ with maximum degree $\nu$.  For example, 
\[
c_{\alpha i}^\ast = c_{\alpha i} + 2 n^{-1/2} c_{\alpha i d} x_d
+ \cdots + (\nu +1) n^{-\nu/2} c_{\alpha i d \cdots d}x_d^\nu.
\]
We do not need to know the coefficients, i.e.  $c_{\alpha i}^\ast, \, c_{\alpha i_1 i_2}^\ast,\ldots$,  of these polynomials in explicit form.  All we need to know is that  (i) the total number of coefficients $c_{\alpha i}^\ast, c_{\alpha i_1 i_2}^\ast, \ldots$ of these polynomials remains bounded as $n \rightarrow \infty$;  and (ii) each of these coefficients remain uniformly bounded as $n \rightarrow \infty$.   Property (i) holds because $q$ and $d$ remain fixed as $n \rightarrow \infty$.  Property (ii) holds because, in each case,
\[
c_{\alpha \cdots}^\ast = c_{\alpha \cdots} + O\left (\frac{\log (n)}{n^{1/2}} \right )
\]
and, by assumption, the $c_{\alpha \cdots}$ are uniformly bounded by $b$.  

Using this result, we define 
\[
W^\ast = \left ( \sum_{i=1}^{d-1} c_{\alpha i}^\ast c_{\beta i}^\ast
\right )_{\alpha , \beta =1}^q,
\]
recalling that we have chosen the labeling so that, in the notation of Lemma 1, the $W_j$ with the largest minimum eigenvalue $\lambda_{\textrm{min}}(W_j)$ is $W_d=W_0-v_d v_d^\top$ where $v_d=(c_{1d}, \ldots c_{qd})^\top$.  Hence, by Lemma 1 there exists a $\theta \in (0,2)$ such that, for $n$ sufficiently large,
\[
\lambda_{\min}(W_d) \geq \frac{\theta}{b} \hskip 0.2truein \textrm{and} \hskip 0.2truein \lambda_{\max}(W_d) \leq \frac{b}{\theta}.  
\]
Hence, for $n$ sufficiently large, 
\[
\lambda_{\min}(W^\ast) \geq \frac{\theta}{2b} \hskip 0.2truein \textrm{and} \hskip 0.2truein \lambda_{\max}(W^\ast) \leq \frac{2b}{\theta}.
\]
Making use of the comments earlier in the proof, it is also clear that 
\[
\max_{\alpha=1, \ldots , q} \max_{i_1, \ldots , i_\nu = 1, \ldots d-1} (\vert c_{\alpha i}^\ast \vert, \vert c_{\alpha i_1 i_2}^\ast \vert , \ldots , \vert c_{\alpha i_1, i_2 \ldots ,i_\nu}^\ast \vert ) \leq \frac{2b}{\theta},
\]
given that $\theta \in (0,2)$.
Therefore, following the key step in the proof of the result on page 1446 of Hall (1986), we have
\begin{align*}
&\sup_{B \in \mathcal{B}} P\left \{  \vert \vert Z \vert \vert \leq b\log n \hskip 0.2truein
\textrm{and} \hskip 0.2truein (p_{n1}(Z), \ldots , p_{nq}(Z))^\top \in (\partial B)^{n^{-\beta}}\right \} \\
&\leq \sup_{B \in \mathcal{B}} \, \sup_{\vert x_d \vert \leq b \log n}
P\left \{ \vert \vert Z^\ast \vert \vert \leq b \log n \hskip 0.1truein \textrm{and} \hskip 0.1truein (p_{n1}^\ast(Z^\ast \vert x_d), \ldots , p_{nq}^\ast(Z^\ast \vert x_d))^\top \in (\partial B)^{n^{-\beta}}\right \}\\
& \leq C(\beta, 2b/\theta, d-1, q) n^{-\beta},
\end{align*}
where we have used the requirement that $\mathcal{B}$ is closed under translations in the sense of (\ref{curlyB_1}).
This completes the proof of the proposition.  \hfill $\square$

\section{The sets $\mathcal{E}_1, \ldots , \mathcal{E}_5$ in Theorem 5.1 in Hall (1992) \label{four_lemmas}}

In this section, to avoid having to introduce another index, we will focus initially on a single population and will derive conditions for the lemmas below to hold.  Then, in Theorem 1, which involves $k$ populations, we will assume that all $k$ populations satisfy the conditions derived in this section.

Following the proof of Theorem 5.1 in Hall (1992), the sets $\mathcal{E}_1, \ldots , \mathcal{E}_5$ are defined as follows.  Let $\mathcal{C}=\mathbb{R}^{nd}$ denote the set of samples $\mathcal{X}=(X_1, \ldots , X_n)$ where each $X_i$ is a real $d$-vector.  The sets $\mathcal{E}_1, \ldots , \mathcal{E}_5$ are subsets of $\mathcal{C}$ and are designed so that $1-P(\mathcal{E}_j)=O(n^{-\lambda})$.  The set $\mathcal{E}_1$ is defined as follows (see Hall, 1992, p.245).  For some constant $C_1>0$, which may be chosen to be large,
\begin{equation}
\mathcal{E}_1=\{\mathcal{X} \in \mathcal{C}: \, \xi_\nu \leq C_1, \, \vert \vert \bar{X}- \mu \vert \vert  \leq C_1^{-1}\}
\label{mathcalE_1}
\end{equation} 
where $\xi_\nu $ is the largest coefficient, in absolute value, of the multivariate polynomials $\check{P}_j(x)=\check{P}_j(0,V, \{\hat{\kappa}_r \})(x)$, $1 \leq j \leq \nu$, where the polynomials $\check{P}_j(x)$ are defined in (\ref{check_P_j}).  Note that we have defined the set $\mathcal{E}_1$ slightly differently from Hall (1992).  Our definition avoids the need to consider moments of partial derivatives of the function $A$ defined in (\ref{new_A}); to compensate, we make explicit smoothness assumptions of the function $A$ in (\ref{new_A}).

The sets $\mathcal{E}_2, \mathcal{E}_3$ and $\mathcal{E}_4$ are defined by
\[
\mathcal{E}_2=\left \{\mathcal{X}: E(\vert \vert X^\ast - \bar{X} \vert \vert^{\max(2 \nu +3, d+1)} \vert \mathcal{X}) \leq C_2
\right \},
\]
where $\bar{X}$ is the mean of the original sample (here we are taking $k$, the number of samples, to be $1$), $X^\ast$ is an observation selected randomly from the original sample $\mathcal{X}$, and $C_2$ is a positive constant; for some $C_3>1$,
\[
\mathcal{E}_3=\left \{ \mathcal{X}: \lambda_{\textrm{max}}(\hat{V}_1) \leq C_3 
\hskip 0.2truein \textrm{and} \hskip 0.2truein \lambda_{\textrm{min}}(\hat{V}_1) \geq C_3^{-1} \right \},
\]
where
\begin{equation}
\hat{V}_1=n^{-1} \sum_{i=1}^n (X_i - \bar{X}) (X_i - \bar{X})^\top
\label{samcov}
\end{equation}
is the sample covariance matrix of the original sample and 
$\lambda_1(A) \geq \cdots \geq \lambda_p(A)$ are the ordered eigenvalues of a $p \times p$ symmetric matrix $A$; and
\[
\mathcal{E}_4=\left \{ \mathcal{X}: \lambda_{\textrm{max}}(\hat{V}^{\dagger}) \leq 2
\hskip 0.2truein \mathrm{and} \hskip 0.2truein \lambda_{\textrm{min}}(\hat{V}^\dagger) \geq \frac{1}{2}\right \},
\]
where 
\[
\hat{V}^\dagger = n^{-1} \sum_{i=1}^n X_i^\dagger \left ( X_i^\dagger \right )^\top,
\]
$X_i^\dagger = Y_i -E(Y_i\vert \mathcal{X})$, with
\[
Y_i = \begin{cases}
\hat{V}_1^{-1/2} (X_i^\ast -\bar{X}) & \textrm{if} \hskip 0.2truein
\vert \vert \hat{V}_1^{-1/2}(X_i^\ast - \bar{X}) \vert \vert \leq n^{1/2}\\
0 & \textrm{otherwise}
    \end{cases}.
\]
so that $E(X_i^\dagger \vert \mathcal{X})$ is the zero vector.

For fixed $u \in (0,\infty)$, the set $\mathcal{E}_5$ is defined by
\begin{equation}
\mathcal{E}_5=\left \{\mathcal{X}: \int \vert \chi^\ast(n^{-(1/2)+u}t)- \chi(n^{-(1/2)+u} t)\vert^m \exp(-C \vert \vert t \vert \vert^{1/2})dt \leq n^{-d(u+1)-(\nu+3)/2} \right \}.
\label{setE5}
\end{equation}

Our first lemma is the following.

\noindent \textbf{Lemma 2}.  \textit{For $\lambda \geq 1$ and $\nu \geq 1$, 
\[
1-P(\mathcal{E}_1)=O(n^{-\lambda}),
\]
provided $E[\vert \vert X \vert \vert^{2(\nu+2)\lambda}]<\infty$.}

\noindent \textbf{Proof of Lemma 2}.  The polynomial $\check{P}_j(x)$ depends on cumulants, and therefore moments, up to and including those of order $\nu+2$, but not larger (this claim follows directly from Lemma 7.1 and Lemma 7.2 of Bhattacharya and Rao, 1976).  Using Chebychev's inequality and the moment condition, the result follows, noting that condition $\vert \vert \bar{X} - \mu \vert \vert \leq C_1^{-1}$ requires a weaker moment condition than that stated in the lemma provided $\nu \geq 1$.  \hfill $\square$

Let us now consider $\mathcal{E}_2$.  The relevant result here is the following.  The definition of $r$ in the statement of Lemma 3 below comes from the definition of $\mathcal{E}_2$ at the bottom of page 245 of Hall (1992).

\noindent \textbf{Lemma 3}.  \textit{For $r=\max(2 \nu+1, d+1)$, where $\nu$ appears in the statement of Theorem 1 and $d$ is the dimension of the vectors $X_j$ that appear in the statement of Theorem 1,
\[
1 - P(\mathcal{E}_2) =O(n^{-\lambda}),
\]
provided that $E(\vert \vert X \vert \vert^{2\lambda r})<\infty$}.

To some extent we follow the approach of Hall (1992).  However, we have been unable to follow all the steps in his argument in the final paragraph of page 245.  Below, we provide a more direct and elementary proof whose conclusion seems to be weaker than that implied by Hall's (1992) conclusion near the bottom of page 245. Hall's conclusion seems to imply that Lemma 3 follows if $E(\vert \vert X \vert \vert^{\lambda r})<\infty$, which is weaker than our condition $E(\vert \vert X \vert \vert^{2\lambda r})<\infty$; we find this discrepancy puzzling but prefer to make use of our result, which we are confident is correct.  It seems that either Hall has used correct reasoning we have been  unable to follow or there is a local error in this part of his argument.

\noindent \textbf{Proof of Lemma 3}.  First, we note that
\begin{align*}
E(\vert \vert X^\ast - \bar{X}\vert \vert^r \vert \mathcal{X})
&= n^{-1} \sum_{i=1}^n \left \{\sum_{j=1}^d (X_i^{(j)} - \bar{X}^{(i)})^2 \right \}^{r/2}\\
&\leq 2^r n^{-1} \sum_{i=1}^n \sum_{j=1}^d (\vert X_i^{(j)} \vert^r +\vert \bar{X}^{(j)}\vert^r).
\end{align*}
Therefore, making use of standard results such as
\[
P\left ( \sum_{t=1}^T  X_t >T\right ) \leq \sum_{t=1}^T P(X_t>1)
\]
and using Chebychev's inequality, we have
\begin{align*}
&P\left \{ E(\vert \vert X^\ast - \bar{X} \vert \vert^r \vert \mathcal{X})> 2^r \sum_{j=1}^d \left \{E(\vert X^{j}\vert^r
+ E\vert \bar{X}^{(j)} \vert^r\right \} + 2^{r+1}d \right \}\\
&\leq P\left \{2^r  n^{-1} \sum_{i=1}^n \sum_{j=1}^d (\vert X_i^{(j)} \vert^r +\vert \bar{X}^{(j)}\vert^r)>2^r\sum_{j=1}^d \left \{E(\vert X^{j}\vert^r
+ E\vert \bar{X}^{(j)} \vert^r\right \} + 2^{r+1} d \right \}
\\
&=P\left \{ n^{-1} \sum_{j=1}^d \left \{ \sum_{i=1}^n \left (\vert X_i^{(j)}\vert^r - E(\vert \bar{X}^{(j)}\vert^r \right ) \right \} + \sum_{j=1}^d ( \vert \bar{X}^{(j)} \vert^r - E(\vert \bar{X}^{(j)} \vert^r)   > 2d \right \} \\
& \leq \sum_{j=1}^d \left \{ P\left [ \bigg \vert n^{-1} \sum_{i=1}^n (\vert X_i^{(j)} \vert^r -E(\vert X^{(j)}\vert^r)\bigg \vert >1\right ] + P\left [ \bigg  \vert \vert \bar{X}^{(j)} \vert^r - E(\vert \bar{X}^{(j)} \vert^r ) \bigg \vert >1 \right ]\right \}\\
&\leq \sum_{j=1}^d \left \{E\bigg \vert n^{-1} \sum_{i=1}^n
\vert X_i^{(j)} \vert^r  - E(\vert X^{(j)}\vert^r )\bigg \vert^{2\lambda }   + E\bigg \vert \vert\bar{X}^{(j)} \vert^r - E(\vert \bar{X}^{(j)} \vert^r ) \bigg \vert^{2 \lambda} \right \}\\
&=O(n^{-\lambda}),
\end{align*}
providing the expectations in the penultimate line are finite, which is the case if and only if $E(\vert \vert X \vert \vert^{2 \lambda r})<\infty$.  \hfill $\square$

For $\mathcal{E}_3$ and $\mathcal{E}_4$ we have the following result.

\noindent \textbf{Lemma 4}.  \textit{Suppose that $E(\vert \vert X \vert \vert^{4 \lambda}) < \infty$ for some $\lambda >1$.  Then}
\[
1-P(\mathcal{E}_3)=O(n^{-\lambda}) \hskip 0.2truein \textrm{and}
\hskip 0.2truein 1 - P(\mathcal{E}_4)=O(n^{-\lambda}).
\]

\noindent \textbf{Proof of Lemma 4}.  We first consider the set
$\mathcal{E}_3$.  Let $\hat{V}$ denote the sample covariance matrix defined in (\ref{samcov}) and let $V_0$ denote the corresponding population covariance matrix.  Suppose that
\[
\lambda_1(V_0) \leq \frac{C_3}{2} \hskip 0.2truein
\textrm{and} \hskip 0.2truein \lambda_q(V_0) \geq \frac{2}{C_3}.
\]
Then we may choose $\epsilon>0$ sufficiently small so that
\[
\{ \mathcal{X}: \lambda_1(\hat{V}_1) \geq C_3\} \bigcup \{ \mathcal{X}: \lambda_q(\hat{V}_1) \leq C_3^{-1} \} \subseteq \{ \mathcal{X}: \vert \vert \hat{V}_1 - V_0 \vert \vert >\epsilon \}.
\]
Then, for such a choice of $\epsilon$, it follows from Chebychev's inequality that
\begin{align*}
P(\vert \vert \hat{V} - V_0 \vert \vert >\epsilon)
& = P(n^{1/2} \vert \vert \hat{V} - V_0 \vert \vert > n^{1/2} \epsilon)\\
&\leq \frac{E(n^{1/2}\vert \vert \hat{V} - V_0 \vert \vert^{2\lambda})}{n^{2\lambda/2}\epsilon^{2 \lambda}}\\
&=O(n^{-\lambda}).
\end{align*}
Finally, note that a necessary and sufficient condition for the $2 \lambda$ moment of $\hat{V}$ to be finite is that $E(\vert \vert X\vert \vert^{4 \lambda})<\infty$.  This proves the first statement in Lemma 4.  The proof of the second statement follows using a similar argument applied to  $X_j^\dagger$.  \hfill $\square$

To deal with $\mathcal{E}_5$, moment conditions are not required; all we need is that Cram\'er's condition (\ref{cramer}) is satisfied.  Lemma 5 below follows from Lemma 5.2 in Hall (1992) plus points made in
the discussion on  p. 252 of Hall (1992).  There is effectively  no change in the version of Hall's (1992) Lemma 5.2  given in our Lemma 5, but for completeness and convenience  we state and prove the result here.

\noindent \textbf{Lemma 5}.  \textit{Suppose that $X$ satisfies  Cram\'er's condition (\ref{cramer}).  Given $d$, and for any $u>0$ and $\nu>0$, if  $m$ is chosen sufficiently large for  the inequality 
\begin{equation}
m>2d(u+1)+\nu+3 + \lambda
\label{mequals}
\end{equation}
to hold, then
\[
1 - P(\mathcal{E}_5) = O(n^{-\lambda}).
\]
}

\noindent \textbf{Proof of Lemma 5}.  As noted on page 252 of Hall (1992), uniformly for $t \in \mathbb{R}^d$,
\begin{equation}
  E[\vert \chi^\ast(t) - \chi(t)\vert^m]  \leq C_1n^{-m/2},
  \label{cfbound}
\end{equation}
for some constant $C_1$ independent of $n$.
Hence, applying Chebychev's inequality to the complement of the set $\mathcal{E}_5$ defined in
(\ref{setE5}), and choosing $m$ sufficiently large to satisfy (\ref{mequals}), it is seen that
\begin{align*}
1-P(\mathcal{E}_5) &= P\bigg \{\int_{t \in \mathbb{R}^d} \vert \chi^\ast(n^{-(1/2)+u}t) - \chi(n^{-(1/2)+u} t)\vert^m
\exp(-C_2 \vert \vert t \vert \vert^{1/2}) dt \\
&\hskip 1.5truein >n^{-d(u+1) - (\nu+3)/2} \bigg \}\\\\
&\leq \frac{E \left \{\int_{t \in \mathbb{R}^d}  \vert \chi^\ast(n^{-(1/2)+u}t) - \chi(n^{-(1/2)+u} t)\vert^m
\exp(-C_2 \vert \vert t \vert \vert^{1/2}) dt  \right \}}
{n^{-d(u+1-(\nu+3)/2}}\\\\
&=\frac{\int_{t \in \mathbb{R}^d} \left \{E\vert \chi^\ast(n^{-(1/2)+u}t) - \chi(n^{-(1/2)+u}t) \vert^m \right \} \exp(-C_2\vert \vert t \vert \vert^{1/2})dt}{n^{-d(u+1) - (\nu+3)/2}}\\\\
&=O \left ( \frac{n^{-m/2}}{n^{-d(u+1) - (\nu+3)/2}} \right )\\
&=O\left ( \frac{1}{n^{(m/2)-d(u+1) -(\nu+3)/2}} \right )\\
&=O(n^{-\lambda}),
\end{align*}
where $C_2$ is a positive constant independent of $t$ and $n$, $m$ is chosen to satisfy (\ref{mequals}) and we have used (\ref{cfbound}) to move from the third line to the fourth line above.  \hfill $\square$

\section{Proof of Theorem 1.  \label{Proof_of_S1}}  

We are now in a position to give our proof of Theorem 1.  Our proof has two components: the first component, Step 1, is modeled on the proof of Theorem 5.1 in Hall (1992) while the second component, Step 2, is the derivation of a sufficient moment condition for Theorem 1 to hold. The moment condition we derive is based on the results in Section \ref{four_lemmas}.  As noted below, this moment condition is very likely to be a long way from being the best possible.

\subsection{Proof of Step 1.}  

Theorem 5.1 in Hall (1992) is split into 4 parts, part (i) to part (iv).  Part (i) of Hall's (1992) proof is concerned with introducing notation and definitions; see Section \ref{modsABC}.  Part (iii) is concerned with identification of the relevant expansion; the structure of the arguments given by Hall (1992) in part (iii) apply without change to our situation when the new definitions given in Section \ref{modsABC} are employed.  Part (iv) of Hall's (1992) proof of Theorem 5.1 is concerned with the proof of Lemma 5.2 in Hall (1992); this also goes through without change in our more general setting.  This leaves part (ii), which is concerned with deriving the key bounds for the relevant Edgeworth expansion.  The remainder of Step 1 is concerned with proving the analogue of part (ii) in our more general setting, using the definitions given in Section \ref{modsABC}.

In the modification of the proof the most significant change is from the class of sets $\mathcal{R}(t)$, $t \in \mathbb{R}$, defined in (5.12) of Hall (1992), to $\mathcal{R}(B)$, defined in Section \ref{curlyR_B}.  

In part (ii) of the proof of Theorem 5.1 in Hall (1992), the arguments from pages 244-249 are virtually identical in our case, the only difference being that there are $k$ samples as opposed to one sample, and in the $k$-sample case the calculations and inequalities need to be applied to each of the $k$ samples in turn.

The first appearances of $(\partial \mathcal{R}(B))^{2 \delta}$ and $(\partial \mathcal{R}^\dagger (B))^{2 \delta}$ are in (5.19) and (5.20) of Hall (1992).  The next appearance of $(\partial \mathcal{R}^\dagger (B))^{2 \delta}$ is in the definition of $b_\delta$ near the bottom of p. 252 of Hall (1992), namely
\begin{equation}
b_\delta = \sum_{j=0}^{\nu + d} n^{-j/2} \int_{(\partial \mathcal{R}^\dagger (B))^{2 \delta}}\vert P_j(-\phi_{0,\hat{V}^\dagger}: \{\hat{\kappa}_{jr}, j=1, \ldots , k\})
(x) \vert dx.
\label{b_delta}
\end{equation}
Moreover, using the fact that the polynomial $P_j(-\phi_{0,\hat{V}^\dagger}: \{\kappa_{jr}\})$ is a polynomial of degree $3j$ multiplied by $\phi_{0,\hat{V}^\dagger}(x)$, as noted by Hall (1992) at the top of page 253, we have the bound
\begin{equation}
b_\delta \leq C_2\int_{(\partial \mathcal{R}^\dagger(B))^{2\delta} }\left (1+\vert \vert x \vert \vert^{3( \nu + kd)} \right )\phi_{0,\hat{V}^\dagger}(x)dx,
\label{b_delta_1}
\end{equation}
where $x \in \mathbb{R}^{kd}$ and $\hat{V}^\dagger$ is defined in (\ref{hatVhatV^dagger}).  Defining 
\begin{equation}
a=-n^{1/2} \hat{V}^{1/2}E(Y\vert \mathcal{X}_1, \ldots , \mathcal{X}_k),
\label{vector_a}
\end{equation}
where $\hat{V}$ is defined in (\ref{hatVhatV^dagger}), $Y=(Y_1^\top, \ldots , Y_k^\top)^\top \in \mathbb{R}^{kd}$ with each $Y_j$ defined via truncation as in (\ref{eq:Yji}), and
$\mathcal{X}_1, \ldots , \mathcal{X}_k$ are the $k$ original samples.

Using exactly the same argument as that given on page 253 of Hall (1992), it is seen that
\[
\left \{ \hat{V}^{1/2} x: x \in (\partial \mathcal{R}(B))^{2 \delta} \right \} \subseteq (\partial (\mathcal{R}(B)+a))^\eta,
\]
where, for some positive constant $C_3$, $\delta = C_3 \eta$ and
$\mathcal{R}(B)+a$ is the translation of $\mathcal{R}(b)$ given by
\[
\mathcal{R}(B) +a = \{y+a \in \mathbb{R}^{kd}: y \in \mathcal{R}(B)\}.
\]
Using exactly the same argument as given by Hall (1992) given in the middle of page 253, it is seen that $\vert \vert a \vert \vert$
is bounded if $\mathcal{X}_j \in \mathcal{E}^{(j)}$, $j=1, \ldots , k$, from which it follows that, for some $C_4>0$ and $C_5>0$, 
\begin{equation}
b_\delta \leq C_4 \int_{(\partial (\mathcal{R}(B)+a ))^\eta}
\exp \left ( - C_5 \vert \vert x \vert \vert^2 \right )dx,
\label{b_delta_2}
\end{equation}
where $b_\delta$ is defined in (\ref{b_delta_1}).  

Note that we may write
\[
\mathcal{R}(B) + a =n^{1/2}\hat{A}(\bar{X}_1+n^{-1/2}(x_1-a_1), \ldots, \bar{X}_k + n^{-1/2}(x_k-a_k)) \in B,
\]
where $x=(x_1^\top, \ldots , x_k^\top)^\top$ and $a=(a_1^\top , \ldots , a_k^\top )^\top$. Following formula (5.27) of Hall (1992) but with minor differences, we now consider the expansion 
\[
n^{1/2}\hat{A}(\bar{X}+n^{-1/2}(x_1 -a_1), \ldots , \bar{X}_k+n^{-1/2}(x_k-a_k)=
\hat{A}_1(x) + \hat{A}_2(x),
\]
where $\hat{A}_u(x)=(A_{u1}(x_1), \ldots ,  A_{uk}(x_k))^\top$, $u=1,2$; and, for $\alpha=1, \ldots , k$,
\begin{align*}
\hat{A}_{1 \alpha}&=\sum_{j_1=1}^k \sum_{i_1=1}^d \hat{A}_{\alpha j_1 i_1}x_{j_1 i_1} 
+\frac{1}{2!} n^{-1/2}
\sum_{j_1, j_2}^k \sum_{i_1, i_2=1}^d \hat{A}_{\alpha j_1 j_2 i_1 i_2}x_{j_1 i_1} x_{j_2 i_2} + \cdots \\
& \hskip 0.3truein  + \frac{1}{(\nu+2)!} n^{-(\nu+1)/2} \sum_{j_1, \cdots , j_{\nu+2}=1}^k 
\, \sum_{i_1, \cdots  , i_{\nu+2}=1}^d 
\hat{A}_{\alpha j_1 \cdots j_{\nu+2} i_1 \cdots i_{\nu+2}} x_{j_1 i_1}  \cdots x_{j_{\nu+2} i_{\nu+2}},
\end{align*}
for $1 \leq s \leq \nu+3$, 
\[
\hat{A}_{1\alpha j_1 \cdots j_s i_i \cdots i_s}=\frac{\partial^s}{\partial x_{j_1 i_1} \cdots \partial x_{j_s i_s}} n^{1/2}\hat{A}(\bar{X}_1+n^{-1/2} x_1, \ldots , \bar{X}_k+n^{-1/2}x_k)
\]
and, uniformly on the set $\mathcal{E} \cap \{\vert \vert x \vert \vert \leq \log n\}$,
\[
\vert \vert \hat{A}_2(x) \vert \vert \leq C^\ast n^{-(\nu+2)/2} (\log n)^{\nu+3}.
\]

We now apply Proposition 1 to the multivariate polynomial $\hat{A}_{1\alpha}$, $\alpha =1, \ldots , q$ and this leads to a bound
\[
b_\delta \leq C_1^\ast n^{-(\nu+1)/2};
\]
see the definition (\ref{b_delta_1}) and the bound (\ref{b_delta_2}).  The remainder of the proof of our Step 1  is similar to the remainder of part (ii) of the proof of Theorem 5.1 in Hall (1992); see in particular page 255 of Hall (1992), recalling that parts (iii) and (iv) of Hall's (1992) proof apply without change in our setting.

\subsection{Proof of Part 2.} 
From Lemma 2, Lemma 3 and Lemma 4 in Section \ref{four_lemmas}, the number of moments required is the largest of
\[
2(\nu+2) \lambda, \hskip 0.2truein 2 \lambda \max(2\nu+3, d+1) \hskip 0.2truein \hbox{and} \hskip 0.2truein 4 \lambda,
\]
which is clearly $2 \max(2 \nu+3, d+1)\lambda$, since $\lambda, \nu, d \geq 1$.  Note that Lemma 5 is not relevant here, as Lemma 5 depends on sufficient smoothness being present (i.e. Cram\'er's condition holding for the underlying populations) but not on the existence of moments.  \hfill $\square$

\section{Further discussion}

We conclude with some discussion.

\noindent 1.  The function $A$ in the statement of Theorem 1 is now a function of $\bar{X}_1^\ast, \ldots , \bar{X}_k^\ast$ whose range is contained in an open subset of $\mathbb{R}^q$, where $a\geq 1$.  Hall (1992, Theorem 5.1) covers the case $q=1$.

\noindent 2.  When $q>1$, the sets $\mathcal{R}(t)$, $t \in \mathbb{R}$ defined in Hall (1992, formula 5.12)) needs to be generalized to $B \in \mathcal{B}$ where $\mathcal{B}$ is defined suitably. In many applications of Theorem 1, it is convenient to take $\mathcal{B}$ to be a class of balls in $\mathbb{R}^q$.

\noindent 3.  Following on from Point 2, we need a multivariate generalization of Lemma 5.3 in Hall (1992).  A suitable generalization, which so far as we are aware is new, is given in Proposition 1 above.

\noindent 4.  Hall (1992) did not attempt to determine explicit moment conditions for Theorem 5.1 to hold.  Indeed, Hall (1992, p.240) states that the number of moments required for Theorem 5.1 to hold ``is not specified in Theorem 5.1 since it is most unlikely that the [lower] bound [for the number of moments needed] produced by our proof is anywhere near the best possible''.  Even in the classical setting (i.e. when Cram\'er's condition is satisfied), there are challenging subtleties which arise when deriving minimal moment conditions for Edgeworth expansions to be valid in the smooth function model; see Hall (1987) and Bhattacharya and Ghosh (1988). So far as we are aware, there has been no further progress on deriving minimal moment conditions since these two papers were published nearly 40 years ago.  Consequently,  deriving minimal moment conditions remains an open problem.

\noindent 5.  Theorem 1 is potentially useful whenever a statistic has a $\chi_q^2$ limit distribution.  Appendix B of Fisher et al. (1996) presents an approach to the asymptotic theory for the bootstrap in such settings. Step B.4 in that paper may be justified by Theorem 1.

\vskip 0.2truein

\section*{Acknowledgements}
This work was supported by Australian Research Council grant DP220102232.  I am grateful to Kassel Hingee and Janice Scealy, co-authors on Hingee et al. (2026), for reading and commenting on an earlier draft of this paper.  I would also like to take this opportunity to express my gratitude to Peter Hall for teaching me about the bootstrap (and other topics) many years ago and for his exceptional  kindness and support at an early stage in my career. 



\end{document}